\documentclass[12pt,reqno,a4paper]{amsart}

\usepackage{pstricks}

\usepackage{pst-node}   
\usepackage{pstcol}     
\usepackage{pst-grad}   
\usepackage{pst-char}   
\usepackage{pst-plot}   

\usepackage{graphicx}
\usepackage{color}

\newcommand{\showgrid}{}

\newcommand{\gridon}{\renewcommand{\showgrid}{\psset{subgriddiv=1,griddots=10,gridlabels=6pt}\psgrid}}

\gridon 

\newgray{hellgrau}{0.9}
\newrgbcolor{dummycolor}{0.75 0.85 0.0}

\long\def\comment#1{%
\psshadowbox[%
fillstyle=solid,fillcolor=dummycolor,linewidth=0.05%
]{\hbox to 12.70cm {\vbox{\strut {\bf Kommentar:} #1\hfill}}}%
}

\advance\voffset    by -1  cm
\advance\hoffset    by -1.5cm
\advance\textwidth  by  3  cm
\advance\textheight by  1  cm

\sffamily

\def\figref#1{Abbildung~\ref{#1}}
\def\heute{\number\day.~\ifcase\month\or
 J\"anner\or Februar\or M\"arz\or April\or Mai\or Juni\or
 Juli\or August\or September\or Oktober\or November\or Dezember\fi
 \space\number\year}

\newif\ifenglish
\englishtrue

\newtheorem{thm}{\ifenglish Theorem\else Satz\fi}
\newtheorem{pro}{Proposition}

\newtheorem{cor}{\ifenglish Corollary\else Korollar\fi}

\newif\ifslides
\slidesfalse

\def\bltxfig#1#2{
\ifslides
\begin{center}
	{\sl #1\/}
	\bigbreak
\else
\begin{figure}
\caption{#1}\label{#2}
\fi
}
\def\eltxfig{
\ifslides
\end{center}
\else
\end{figure}
\fi
}

\def\bltxtab#1#2{
\ifslides
\begin{center}
	{\sl #1\/}
	\bigbreak
\else
\begin{table}
\caption{#1}\label{#2}
\fi
}
\def\eltxtab{
\ifslides
\end{center}
\else
\end{table}
\fi
}

\def\defeq{:=}
\def\absof#1{\left|#1\right|}
\def\setof#1{\left\{#1\right\}}
\def\of#1{\!\left(#1\right)}

\def\pas#1{\left(#1\right)}
\def\brk#1{\left[#1\right]}

\def\bit{\begin{itemize}}
\def\eit{\end{itemize}}
\def\beq{\begin{equation}}
\def\eeq{\end{equation}}

\def\C{{\mathbb C}}

\def\Z{{\mathbb Z}}

\def\i{{\mathbf i}}

\def\e{{\mathbf e}}
\def\1{{\mathbf 1}}

\def\dd{{\rm d}}

{\begin{list}{#1}{\parsep0em \itemsep0.3em \labelwidth1em
\labelsep0.5em \leftmargin1.5em }} {\end{list}}

\def\lgt{Lindstr\"om--Gessel--Viennot}

\def\ddd{m}	
\def\OO{O}	
\def\oo{o}	
\def\ndivs#1{{\mathrm d}\of{#1}}	
\def\pochhammer#1#2{\pas{#1}_{#2}}
\def\fallfac#1#2{\left(#1\right)_{\left.#2\right\downarrow}}

\def\mytheta{{\bar\vartheta}}

\def\singlesum{{\mathbf S_1}}
\def\doublesum{{\mathbf S_2}}
\def\dee{{\mathrm d}}
\def\deeof#1{{\mathrm d} #1}
\def\differ#1#2#3{\frac{\dee^{#3}}{\dee #2^{#3}}#1}
\def\Iverson#1{\brk{#1}}

\def\path{{%\mathrm 
P}}
\def\figref#1{Figure~\ref{#1}}

\def\minheight{\underline{m}}
\def\bigpas#1{\Bigl(#1\Bigr)}
\def\bdetof#1{\Bigl|#1\Bigr|}

\def\pochhammer#1#2{\pas{#1}^{\overline{#2}}}
\def\fallfac#1#2{\pas{#1}^{\underline{#2}}}

\usepackage{epsfig} 

\begin{document}

\bibliographystyle{plain}

\title{Asymptotics of the average height of $2$--watermelons with a wall.}

\begin{abstract}
We generalize the classical work of de~Bruijn, Knuth and Rice
(giving the asymptotics of the average height of Dyck paths of
length $n$) to the case of $p$--watermelons with a wall (i.e., to a
certain family of $p$ nonintersecting Dyck paths; simple
Dyck paths being the special case $p=1$.) We work out this
asymptotics for the case $p=2$ only, since the computations
involved are already quite complicated (but might be of some interest
in their own right).
\end{abstract}

\author{Markus Fulmek}
\address{Fakult\"at f\"ur Mathematik, 
Nordbergstra\ss e 15, A-1090 Wien, Austria}
\email{{\tt Markus.Fulmek@Univie.Ac.At}\newline\leavevmode\indent
{\it WWW}: {\tt http://www.mat.univie.ac.at/\~{}mfulmek}
}

\date{\today}
\thanks{
Research supported by the National Research Network ``Analytic
Combinatorics and Probabilistic Number Theory'', funded by the
Austrian Science Foundation. 
}

\maketitle

\section{Introduction}
The model of vicious walkers was originally introduced by Fisher
\cite{fisher:walks} and received much interest, since it leads
to challenging enumerative questions. Here, we consider special
configurations of vicious walkers called {\em $p$--watermelons with
a wall\/}.

Briefly stated, a $p$--watermelon of length $n$ is a family
$\path_1,\dots \path_p$ of $p$ nonintersecting lattice paths in
$\Z^2$, where
\bit
\item $\path_i$ starts at $\pas{0,2i-2}$ and ends at $\pas{2n,2i-2}$,
	for $i=1,\dots,p$,
\item all the steps are directed north--east or south--east, i.e.,
	lead from lattice point $\pas{i,j}$ to $\pas{i+1,j+1}$ or to
	$\pas{i+1,j-1}$,
\item no two paths $\path_i$, $\path_j$ have a point in common (this
	is the meaning of ``nonintersecting'').
\eit

The {\em height\/} of a $p$--watermelon is the $y$--coordinate
of the highest lattice point contained in any of its paths (since the paths are nonintersecting, it
suffices to consider the lattice points contained in the highest path
$\path_p$; see \figref{fig:watermelon} for an illustration.)

\begin{figure}
\caption{A $6$--watermelon of length 46 and height 20}
\label{fig:watermelon}
\bigskip
\psset{unit=0.2cm}
\pspicture(0,-10)(45,22)
%
%
\psline[linewidth=0.1](0,0)(46,0)
\psline[linewidth=0.1](0,-10)(0,20)
\psline[linewidth=0.05](0,0)(2,-2)(3,-1)(8,-6)(10,-4)(16,-10)(19,-7)(21,-9)(26,-4)(31,-9)(38,-2)(41,-5)(46,0)
\psline[linewidth=0.05](0,2)(2,0)(3,1)(6,-2)(8,0)(16,-8)(20,-4)(22,-6)(29,1)(34,-4)(40,2)(42,0)(44,2)(45,1)(46,2)
\psline[linewidth=0.05](0,4)(2,2)(3,3)(6,0)(9,3)(11,1)(12,2)(13,1)(15,3)(22,-4)(32,6)(34,4)(35,5)(37,3)(38,4)(39,3)(41,5)(43,3)(45,5)(46,4)
\psline[linewidth=0.05](0,6)(2,4)(4,6)(5,5)(7,7)(9,5)(16,12)(19,9)(20,10)(21,9)(27,3)(32,8)(34,6)(40,12)(46,6)
\psline[linewidth=0.05](0,8)(2,6)(4,8)(5,7)(7,9)(8,8)(15,15)(16,14)(17,15)(18,14)(19,15)(25,9)(30,14)(32,12)(33,13)(35,11)(36,12)(38,14)(39,13)(40,14)(46,8)
\psline[linewidth=0.05](0,10)(10,20)(12,18)(13,19)(15,17)(17,19)(19,17)(22,20)(28,14)(32,18)(31,17)(33,19)(34,18)(36,20)(46,10)
\pscircle(0,0){0.25}
\pscircle(0,2){0.25}
\pscircle(0,4){0.25}
\pscircle(0,6){0.25}
\pscircle(0,8){0.25}
\pscircle(0,10){0.25}
\pscircle(46,0){0.25}
\pscircle(46,2){0.25}
\pscircle(46,4){0.25}
\pscircle(46,6){0.25}
\pscircle(46,8){0.25}
\pscircle(46,10){0.25}
\pscircle[linecolor=red,fillstyle=solid,fillcolor=red](10,20){0.25}
\pscircle[linecolor=red,fillstyle=solid,fillcolor=red](22,20){0.25}
\pscircle[linecolor=red,fillstyle=solid,fillcolor=red](36,20){0.25}
\psline[linecolor=red,linewidth=0.025](5,20)(41,20)
\rput[l](42.5,20){{\red\small 
	level $20$}} 
\endpspicture
\end{figure}

A $p$--watermelon of length $n$ {\em with a wall\/} has the additional
property that none of the paths ever goes below the line $y=0$ (since
the paths are nonintersecting, it suffices to impose this condition on
the lowest path $\path_1$; see \figref{fig:watermelon-wall} for an illustration.).

\begin{figure}
\caption{A $3$--watermelon with a wall of length 11 and height 12}
\label{fig:watermelon-wall}
\bigskip
\psset{unit=0.3cm}
\pspicture(0,0)(22,14)
%
%
\psline[linewidth=0.05](0,0)(22,0)
\psline[linewidth=0.05](0,0)(0,14)
\psline[linewidth=0.05](22,0)(22,14)
\psline[linewidth=0.05](0,14)(22,14)
\psline[linewidth=0.025](0,0)(2,2)(3,1)(5,3)(6,2)(10,6)(16,0)(18,2)(20,0)(21,1)(22,0)
\psline[linewidth=0.025](0,2)(3,5)(4,4)(10,10)(15,5)(16,6)(18,4)(19,3)(20,4)(22,2)
\psline[linewidth=0.025](0,4)(4,8)(5,7)(10,12)(15,7)(16,8)(19,5)(20,6)(21,5)(22,4)
\pscircle(0,0){0.25}
\pscircle(0,2){0.25}
\pscircle(0,4){0.25}
\pscircle(22,0){0.25}
\pscircle(22,2){0.25}
\pscircle(22,4){0.25}
\pscircle[linecolor=red,fillstyle=solid,fillcolor=red](10,12){0.25}
\psline[linecolor=red,linewidth=0.025](7,12)(13,12)
\rput[l](13.5,12){{\red\small 
	level $12$}} 
\endpspicture
\end{figure}

In \cite{bonichon-mosbah:watermelons}, Bonichon and Mosbah 
considered (amongst other things) the {\em average height\/}
$H\of{n,p}$ of $p$--watermelons of length $n$ with a wall,
\begin{multline*}
H\of{n,p} = \frac{1}{\#\of{\text{all $p$-watermelons of length $n$}}}
\\
\times\sum_h h\cdot {
	\#\of{\text{all $p$-watermelons of length $n$ and height $h$}}
	},
\end{multline*}
and derived by computer experiments the
following conjectural asymptotics \cite[4.1]{bonichon-mosbah:watermelons}:
\begin{equation}
\label{eq:bonichon-mosbah}
H\of{n,p} \sim \sqrt{\pas{1.67 p-0.06} 2 n} + \oo\of{\sqrt n}.
\end{equation}

The purpose of this paper is to work out the exact asymptotics for the
simple special case $p=2$. This will be done by imitating the
classical reasoning of de Bruijn, Knuth and Rice \cite{deBruijn+Knuth+Rice} for the case $p=1$ (i.e., for the average height of Dyck
paths). However, even the case $p=2$ involves rather complicated
computations. In particular, we shall need informations about
residues and evaluations of a double Dirichlet series, which we
shall (partly) obtain by imitating Riemann's representation
of the zeta function
\cite[section 1.12, (16)]{erdelyi:high-spec-func}.

\subsection{Notational conventions}
For $k,n\in\Z$, we shall use the notation introduced
in \cite{GrahamKnuthPatashnik} for the
rising and falling factorial powers, i.e.
\begin{equation*}
\pochhammer{n}{k}\defeq 0 \text{ if } k < 0,\;
\pochhammer{n}{0}\defeq 1,\;
	\pochhammer{n}{k}\defeq n\cdot\pas{n+1}\cdots\pas{n+k-1}
	\text{ if } k > 0
\end{equation*}
and
\begin{equation*}
\fallfac{n}{k}\defeq 0 \text{ if } k < 0,\;
\fallfac{n}{0}\defeq 1,\;
	\fallfac{n}{k}\defeq n\cdot\pas{n-1}\cdots\pas{n-k+1}
	\text{ if } k > 0.
\end{equation*}
For the binomial coefficient we adopt the convention
\begin{equation*}
\binom{n}{k} \defeq
	\begin{cases}
	\frac{\fallfac{n}{k}}{k!} & \text{ if } 0\leq k \leq n,\\
	0 & \text{ else.}
	\end{cases}
\end{equation*}
Moreover, we shall use Iverson's notation:
\begin{equation*}
\Iverson{\text{some assertion}} =
\begin{cases}
	1 &\text{ if ``some assertion'' is true}, \\
	0 &\text{ else.}
\end{cases}
\end{equation*}
\subsection{Organization of the material presented}
This paper is organized as follows:
\bit
\item In Section~\ref{sec:enumeration}, we present
	{\em exact} enumeration formulas for the average
	height of $p$--watermelons with a wall in terms of
	certain determinants. Moreover,
	we make these formulas more explicit
	(in terms of sums of binomial coefficients)
	for the simple cases $p=1$ and $p=2$.
\item In Section~\ref{sec:asymptotics}, we first review the
	classical reasoning for the asymptotics of the average
	height of $1$--watermelons with a wall, which was given
	by de Bruijn, Knuth and Rice \cite{deBruijn+Knuth+Rice}.
	Then we show how this reasoning can be modified for the case of
	$2$--watermelons with a wall.
\item In appendix~\ref{sec:appendix}, 
	we summarize background
	information on
	\bit
	\item Stirling's approximation,
	\item certain residues and values of the gamma and
		zeta function,
	\item a certain double Dirichlet series
		and Jacobi's theta function
	\eit
	which are needed in our presentation.
\eit

\subsection{Acknowledgements}
I am very grateful to Professor Kr\"atzel for pointing out
to me
how the poles and residues of
certain Dirichlet series can be obtained in a simple way by
using the reciprocity law for Jacobi's theta function, and to
Christian Krattenthaler for many helpful discussions.

\section{Exact enumeration}
\label{sec:enumeration}
For a start, we gather some exact enumeration results.

\subsection{The number of $p$--watermelons with a wall}
We have the following generalization of the enumeration of
$1$--watermelons with a wall (i.e., Dyck paths) of length $n$,
which is given by the Catalan numbers
\begin{equation*}
C\of{n}=\frac{1}{n+1}\binom{2n}{n}.
\end{equation*}

\begin{pro}
The number $C\of{n,p}$ of all $p$--watermelons with a wall of
length $n$ is given as
\begin{equation}
\label{eq:all-watermelons-prod}
C\of{n,p}=\prod_{j=0}^{p-1}
\frac{
	\binom{2n+2j}{n}
}{
	\binom{n+2j+1}{n}
}.
\end{equation}
\end{pro}
\begin{proof}
This is a special case of Theorem 6 in \cite{kratt-guttmann-viennot}.
\end{proof}

\subsection{$1$--watermelons with a wall and height restrictions}
In order to obtain the {\em average\/} height, we count $p$--watermelons with a wall of length $n$
which do not {\em exceed\/} height $h$.

To this end, we employ the following formula
(see \cite[p.~6, Theorem~2]{mohanty:lattice-path:1979}):
\begin{thm}
\label{thm:mohanty}
Let $u$, $d$ be nonnegative integers, and let $b$, $t$ be
positive integers, such that $-b<u-d<t$.
The number of lattice paths from $\pas{0,0}$ to $\pas{u+d,u-d}$,
which do not touch neither line $y=-b$ nor line $y=t$, equals
\begin{equation}
\label{eq:mohanty}
\sum_{k\in\Z} \pas{
	\binom{u+d}{u-k\pas{b+t}}
	-
	\binom{u+d}{u-k\pas{b+t}+b}
}.
\end{equation}
\end{thm}
\begin{cor}
Let $i$, $j$ and $h$ be integers such that $0\leq 2i,2j\leq h$.
The number 
of all lattice paths from $\pas{0,2i}$ to
$\pas{2n,2j}$ which lie {\em between\/} the lines $y=0$ and
$y=h\geq 0$ is
\begin{equation}
\label{eq:all-h-gen-dyck-paths}
\ddd\of{n,i,j,h}\defeq
\sum_{k\in\Z}
\pas{
	\binom{2n}{n-i+j-k\pas{h+2}}-
	\binom{2n}{n+i+j-k\pas{h+2}+1}
}.
\end{equation}
The special case $i=j=0$ can be written as
\begin{multline}
\label{eq:all-h-dyck-paths-folded}
\ddd\of{n,h}\defeq \ddd\of{n,0,0,h}=
\frac1{n+1}\binom{2n}{n} - \\
\sum_{k\geq1}
\pas{
	\binom{2n}{n-k\pas{h+2}-1}-
	2 \binom{2n}{n-k\pas{h+2}} +
	\binom{2n}{n-k\pas{h+2}+1}
}.
\end{multline}
\end{cor}
\begin{proof}
Set $u=n-i+j$, $d=n+i-j$, $b=2i+1$ and $t=h+1-2i$ in
\eqref{eq:mohanty}.
\end{proof}
\begin{cor}
For $n>0$,
the number of all $p$--watermelons with a wall of length $2n$, which
do not exceed height $h$, is given by the following determinant:
\begin{equation}
\label{eq:all-h-watermelons}
C\of{n,p,h}=\bdetof{m\of{n,i,j,h}}_{i,j=0}^{p-1}.
\end{equation}
\end{cor}
\begin{proof}
For $h>2p-2$ this follows by a direct application of the
\lgt\ method \cite{gessel-viennot:det}.

For $0\leq h\leq 2p-2$, the determinant equals $0$ (as it should),
since $\ddd\of{n,i,j,2i}=0$ and
$\ddd\of{n,i,j,2i+1}+\ddd\of{n,i+1,j,2i+1}=0$
for all $j$.
\end{proof}


\subsection{The average height of $1$--watermelons with a wall}
The following is a condensed version of the reasoning given
in \cite{deBruijn+Knuth+Rice}. 
Denote by $\minheight\of{n,h}$ the number of all Dyck paths, starting
at $(0,0)$ and ending at $(2n,0)$, which reach {\em at least\/}
height $h$. By \eqref{eq:all-h-dyck-paths-folded}, we obtain
\begin{align*}
\minheight\of{n,h} &= \ddd\of{n,n} - \ddd\of{n,h-1}
 = C\of{n,1} - C\of{n,1,h-1}\\ 
&= \sum_{k\geq1}
\pas{
	\binom{2n}{n-k\pas{h+1}-1}-
	2 \binom{2n}{n-k\pas{h+1}} +
	\binom{2n}{n-k\pas{h+1}+1}
}. 
\end{align*}

So, the average height of $1$--watermelons with a wall (i.e.,Dyck
paths) of length $n$ is
\begin{align}
H\of{n,1}
&= \frac{1}{C_n}\sum_{h=1}^{n} \minheight\of{n,h}
	= \frac{1}{C_n}\pas{-C_n+
		\sum_{h=0}^{n} \minheight\of{n,h}}
	\notag\\
&= 
-1+
\pas{n+1}
\pas{\sum_{k\geq1}
	\ndivs k\frac{
		\pas{
			\binom{2n}{n-k-1}-
			2 \binom{2n}{n-k} +
			\binom{2n}{n-k+1}
		}
	}{
		\binom{2n}{n}
	}
},\label{eq:1-watermelons-average}
\end{align}
where $d\of k$ denotes the number of positive divisors of
$k$. Introducing the notation
\begin{equation}
\label{eq:S-simple}
S\of{n,a}\defeq\sum_{k\geq 1}
	\ndivs k\frac{\binom{2n}{n-k+a}}{\binom{2n}{n}},
\end{equation}
we arrive at
\begin{equation}
\label{eq:avg-1-watermelons}
H\of{n,1} = \pas{n+1}\bigpas{S\of{n,1}-2\,S\of{n,0}+S\of{n,-1}}-1,
\end{equation}
which is equivalent to 
equation (23) in \cite{deBruijn+Knuth+Rice}.

\subsection{The average height of $p$--watermelons with a wall}
In generalization of the above notation, denote by
$\minheight\of{n,p,h}$
the number of all $p$--watermelons of
length $n$, which reach {\em at least\/}
height $h$, i.e.,
\begin{equation*}
\minheight\of{n,p,h} = C\of{n,p} - C\of{n,p,h-1}.
\end{equation*}
Clearly, we have the following exact formula for the average height
of $p$--watermelons of length $n$:
\begin{equation}
\label{eq:p-watermelons-average}
H\of{n,p}=
\frac{1}{C\of{n,p}}\sum_{h=1}^{n+2p-2} \minheight\of{n,p,h}.
\end{equation}

\subsubsection{The average height of $2$--watermelons with a wall}
Let
\begin{equation}
S\of{n,a,b} =
\sum_{j\geq 1}
   \sum_{k\geq 1}
     \dd\of{\gcd\of{j,k}}\,
     \frac{
     	\binom{2\,n}{n - j + a}\,
        \binom{2\,n}{n - k + b}
       }{\binom{2\,n}{n}^2}.
\end{equation}
From \eqref{eq:all-h-watermelons}, straightforward (but rather
tedious) computations lead to the following formula:
\begin{equation}
\label{eq:all-2-watermelons-normal}
H\of{n,2}=
\frac{\pochhammer{n+1}{2}}{12\pas{2n+1}}
\bigpas{
	\pochhammer{n+1}{3}\doublesum\of n +
	\singlesum\of n
}
-1,
\end{equation}
where
\begin{align}
\singlesum\of n
&=
-20 (n-1) (n+2) S(n,0)+15 n (n+1) (S(n,-1)+S(n,1))+ \notag\\
&\phantom{=}
(n+3)(6 S(n,-1)-16 S(n,0)+6 S(n,1))+(n-2) (6 S(n,-1)+ \notag\\
&\phantom{=}
8 S(n,0)+6 S(n,1))-6 n (n+3) (S(n,-2)+S(n,2))+\notag\\
&\phantom{=}
(n+2)(n+3) (S(n,-3)+S(n,3)),\label{eq:singlesum}\\
\doublesum\of n
&=
S(n,-2,-2)-S(n,-1,-3)-2 S(n,-1,-2)+S(n,-1,-1)+2S(n,-1,0)- \notag\\
&\phantom{=}
S(n,-1,3)+2 S(n,0,-3)-4 S(n,0,0)+2S(n,0,3)-S(n,1,-3)-\notag\\
&\phantom{=}
2 S(n,1,-2)+2 S(n,1,-1)+2 S(n,1,0)+S(n,1,1)-S(n,1,3)+\notag\\
&\phantom{=}
2 S(n,2,-2)-2 S(n,2,-1)-2 S(n,2,1)+S(n,2,2).\label{eq:doublesum}
\end{align}
 
\section{Asymptotic enumeration}
\label{sec:asymptotics}
\subsection{Asymptotics of the average height of $1$--watermelons}
\label{sec:asym-avg-1}
In the case of $1$--watermelons, the asymptotic of the average
height \eqref{eq:1-watermelons-average} is well--known, see
\cite[Proposition 7.7]{flajolet-sedgewick:mellin}
or \cite[equation (34)]{deBruijn+Knuth+Rice}:
\begin{equation}
\label{eq:1-watermelons-average-asympt}
H\of{n,1}\simeq\sqrt{\pi n}-\frac32+
\OO\of{n^{-\frac12+\epsilon}}.
\end{equation}
We repeat the classical reasoning of de~Bruijn, Knuth
and Rice \cite{deBruijn+Knuth+Rice}, in order to make clear
the basic idea, which we shall also employ for the case
$p=2$ later.

\noindent{\em Proof of \eqref{eq:1-watermelons-average-asympt}:\/}
From \eqref{eq:avg-1-watermelons} it is clear that we need to
investigate the asymptotic behaviour of
$S\of{n,a}=
\sum_{k=1}^n \ndivs k\frac{\binom{2n}{n-k-a}}{\binom{2n}{n}}$.
Note that the sums $S\of{n,a}$ in \eqref{eq:avg-1-watermelons}
are multiplied with
a factor of order $1$. So if we are interested in the
asymptotics of $H\of{n,1}$ up to some $\OO\of{n^{-\alpha}}$,
we need the asymptotics for $S\of{n,a}$ up to $\OO\of{n^{-\alpha-1}}$;
for our case, $\alpha=1-\epsilon$ is sufficient.

The basis of the following considerations
is the asymptotic expansion of the quotient of binomial
coefficients \eqref{eq:binkoeffquotapprox} (see appendix~\ref{sec:stirlingsformula}).


\subsubsection{The asymptotics of $S\of{n,a}$ for $a$ fixed, $n\to\infty$}
\label{sec:asymSsingle}
First we observe that
\begin{equation*}
\frac{\binom{2n}{n+a-k}}{\binom{2n}{n}}
=
\OO\of{\exp\of{-n^{2\epsilon}}}
\end{equation*}
if $\absof{
\frac{k-a}{n}}\geq n^{\epsilon-\frac12}$, i.e., if
$k\geq n^{\frac12+\epsilon} + a$ (see \eqref{eq:binkoeffquotapprox}
in appendix~\ref{sec:stirlingsformula}). Therefore, the sum of
all terms with $k\geq n^{\frac12+\epsilon} + a$ is negligible in
\eqref{eq:S-simple}, being $\OO\of{n^{-m}}$ for all $m>0$, and we
may take $
\frac{k-a}{n}=\OO\of{n^{\epsilon-\frac12}}$ in
\eqref{eq:binkoeffquotapprox}.



Next, we take \eqref{eq:binkoeffquotapprox} up to order $n^{-1}$
and substitute $x\to\frac{k-a}{n}$: Pulling out the leading term
$\e^{\frac{-k^2}{n}}$ and expanding the rest with respect
to $k$ gives
\begin{multline}
\label{eq:binomial-quotient-asymptotics}
\frac{\binom{2n}{n-k-a}}{\binom{2n}{n}}=
{\e^{\frac{-k^2}{n}}}
\Bigl(
	1 - \frac{a^2}{n} +
    k\,\pas{\frac{2\,a}{n} - \frac{a + 2\,a^3}{n^2}} + 
    \frac{\left( 1 + 4\,a^2 \right) \,k^2}{2\,n^2}  + \\
    \frac{\left( 5\,a + 4\,a^3 \right) \,k^3}{3\,n^3} - 
    \frac{k^4}{6\,n^3}-
    \frac{a\,k^5}{3\,n^4}  \Bigr)
    +
    {\OO\of{n^{-2+
    \epsilon}}}.
\end{multline}
Now we consider the following function
\begin{equation}
g\of{n,b}\defeq\sum_{k\geq 1}k^b\,\ndivs{k}\e^{-k^2/n},
\end{equation}
and observe that here the terms for
$k\geq n^{1/2+\epsilon}$ are {\em again\/} negligible:
\begin{equation*}
\sum_{k\geq n^{1/2+\epsilon}}k^b\,\ndivs{k}\e^{-k^2/n}
=\OO\of{n^{-m}}\text{ for all } m>0.
\end{equation*}
Hence we directly obtain from \eqref{eq:binomial-quotient-asymptotics}:
\begin{multline}
S\of{n,a} = 
\pas{1-\frac{a^2}{n}}g\of{n,0}
+\pas{\frac{2a}{n}-\frac{2a^3+a}{n^2}}g\of{n,1}
+\pas{\frac{4a^2+1}{2n^2}}g\of{n,2}\\
+\pas{\frac{4a^3+5a}{3n^3}}g\of{n,3}
-\pas{\frac{1}{6n^3}}g\of{n,4}
-\frac{a}{3n^4}\,g\of{n,5}
+\OO\of{n^{-2+\epsilon}g\of{n,0}}.
\end{multline}
(This is equation (27) in \cite{deBruijn+Knuth+Rice}.) Note
that the coefficients for $g\of{n,k}$ are odd functions of $a$
for odd $k$ and obtain:
\begin{equation}
\label{eq:avg-1-watermelon-asym1}
S\of{n,1}-2\,S\of{n,0}+S\of{n,-1} 
=-\frac{2}{n}\,g\of{n,0}+\frac{4}{n^2}\,g\of{n,2}+
	\OO\of{n^{-2+\epsilon}g\of{n,0}}.
\end{equation}
So we reduced our problem to that of obtaining an asymptotic
expansion for $g\of{n,b}$. Note that
we need this information only for $b$ even. It follows from the
computations presented in appendix~\ref{sec:gnb-asym}, that
$g\of{n,2b}=\OO\of{n^{b+1/2}\log\of{n}}$, and that we have
for all $m\geq 0$:
\begin{align}
g\of{n,0} &=
	\frac14\sqrt{\pi n}\log\of n +
	\pas{\frac34\gamma -\frac12\log\of 2}\sqrt{\pi n}
	+\frac14 +\OO\of{n^{-m}},
\notag \\
g\of{n,2} &= 
	\frac{n}8\sqrt{\pi n}\log\of n +
	\pas{\frac14+\frac38\gamma-\frac14\log\of 2}n\sqrt{\pi n}+
	\OO\of{n^{-m}}.\label{eq:gn2}
\end{align}
Inserting this information in \eqref{eq:avg-1-watermelon-asym1}
we immediately obtain 
the desired result \eqref{eq:1-watermelons-average-asympt}.
\hfill\qedsymbol

\subsection{Asymptotics of the average height of $2$--watermelons}
\begin{figure}
\caption{Illustration of the quotient
$q\of n = \frac{H\of{n,2}}{2.57758\, \sqrt{n} - 2}$, for $n\leq 1000$.}
\label{fig:convergence}
\centerline{\epsfxsize12cm\epsffile{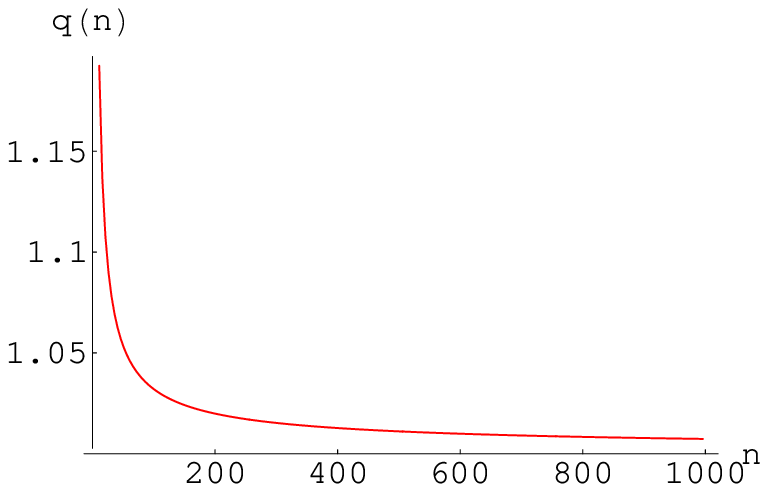}}
\end{figure}

We shall modify the reasoning from section~\ref{sec:asym-avg-1}
appropriately.
%
%
In doing so, it turns out that we have to deal with the double
Dirichlet series
$\sum_{k,l\geq1}\frac{k^{2a}l^{2b}}{\pas{k^2+l^2}^{s}}$ for
integers $a,b\geq 0$. Proposition~\ref{pro:dirichlet-series}
(see section~\ref{sec:dirichlet-series})
states that this series 
is convergent in the half--plane $\Re\of z > a+b+1$ and defines
a meromorphic function $Z\of{a,b;z}:\C\to\C$ which has a simple
pole at $z=a+b+1$, and an additional simple pole at $z=a+b+\frac12$
only if $a=0$ or $b=0$. Hence, we can write
\begin{equation}
\label{eq:Zdefinition}
 Z\of{a,b;z} = \frac{r_{a,b}}{z-a-b-\frac12} + c_{a,b} +
\OO\of{z-a-b-\frac12}.
\end{equation}
Given this ``implicit'' definition of the
numbers $c_{a,b}$ we will show:
\begin{multline}
\label{eq:2-watermelons-average-asympt}
H\of{n,2}\simeq
\sqrt{\pi n}\Bigl(
-2c_{0,0}+8c_{1,0}-9c_{1,1}-9c_{2,0}+15c_{2,1} \\
	+35c_{2,2}+5c_{3,0}-35c_{3,1}\Bigr)-2+
\OO\of{n^{-\frac12+\epsilon}}.
\end{multline}
Using the representations of the constants $c_{a,b}$ by certain
integrals (see \eqref{eq:absolute_terms} in
appendix~\ref{sec:dirichlet-series}), we obtain the following approximative asymptotics by numerical integration (carried out with {\sc Mathematica})
\begin{equation*}
H\of{n,2}\simeq 2.57758\, \sqrt{n} - 2+
\OO\of{n^{-\frac12+\epsilon}},
\end{equation*}
which conforms well to Bonichon's and Mosbah's conjecture
\eqref{eq:bonichon-mosbah} for the case $p=2$, which yields
approximately $2.56125\, \sqrt{n}$. \figref{fig:convergence}
shows the quotient $q\of n = \frac{H\of{n,2}}{2.57758\, \sqrt{n} - 2}$
for small $n$. For example, $q\of{1000} = 1.00734$.

\noindent{\em Proof of \eqref{eq:2-watermelons-average-asympt}:\/}
Note that in \eqref{eq:all-2-watermelons-normal},
the ``single sums'' $S\of{n,a}$ are multiplied with
a rational function in $n$ of order at most $3$, while the
``double sums'' $S\of{n,a,b}$ are multiplied with a factor of
order $4$: So if we are interested in the
asymptotics of $H\of{n,2}$ up to some $\OO\of{n^{-\alpha}}$,
we need the asymptotics for $S\of{n,a}$ up to $\OO\of{n^{-\alpha-3}}$
and for $S\of{n,a,b}$ up to $\OO\of{n^{-\alpha-4}}$;
for our case, $\alpha=1-\epsilon$ is sufficient.


%
\subsubsection{The asymptotics of $S\of{n,a}$
for $a$ fixed, $n\to\infty$} Basically, we repeat the
computations from section \ref{sec:asymSsingle}. The only
difference is that 
we need higher orders now. After some calculations, we
obtain:
\begin{multline}
\label{eq:singlesumIng}
\singlesum\of n=
-\frac{24 \left(4 n^2+20 n+89\right)
   g\of{n,0}}{n^3}+\frac{4 \left(96
   n^2+1065 n+3656\right)
   g\of{n,2}}{n^4}\\
   -\frac{\left(288
   n^2+4060 n+12213\right)
   g\of{n,4}}{n^5}+\frac{8 \left(8
   n^2+107 n+335\right)
   g\of{n,6}}{n^6}\\
   -\frac{(96 n+521)
   g\of{n,8}}{3 n^7}+\frac{10
   g\of{n,10}}{3 n^8}
   +\OO\of{n^{-4+\epsilon}g\of{n,0}}.
\end{multline}
%
%
%
%
Recall that in appendix~\ref{sec:gnb-asym} it is proved that
$g\of{n,2b}=\OO\of{n^{b+1/2}\log\of{n}}$. Moreover, the arguments
in appendix~\ref{sec:gnb-asym} show that we have (in addition to
\eqref{eq:gn2}) for all $m\geq 0$:
\begin{align*}
g\of{n,4} &= 
	\frac{3 n^2}{16}\sqrt{\pi n}\log\of n +
	\pas{\frac12+\frac9{16}\gamma-\frac38\log\of 2}n^2\sqrt{\pi n}+
	\OO\of{n^{-m}},
\\
g\of{n,6} &= 
	\frac{15 n^3}{32}\sqrt{\pi n}\log\of n +
	\pas{\frac{23}{16}+\frac{45}{32}\gamma-\frac{15}{16}\log\of 2}n^3\sqrt{\pi n}+
	\OO\of{n^{-m}},
\\
g\of{n,8} &= 
\frac{105 n^4}{64} \sqrt{\pi n} \log\of{n} +
	\pas{
  \frac{11}{2}
   +\frac{315}{64} \gamma  
    -\frac{105}{32} \log\of{2}
   } n^4 \sqrt{\pi n} +
	\OO\of{n^{-m}},
\\
g\of{n,10} &= 
 \frac{945 n^5}{128} \sqrt{\pi n} \log\of{n}+
 \pas{
   \frac{1689}{64}
   +\frac{2835}{128} \gamma 
   -\frac{945}{64} \log\of{2}
   }n^5 \sqrt{\pi n}+
	\OO\of{n^{-m}}.
\end{align*}
Inserting this information in \eqref{eq:singlesumIng}
we immediately obtain 
the first part of the desired result:
\begin{equation}
\label{eq:h2firstpart}
\frac{\pochhammer{n+1}{2}}{12\pas{2n+1}}\singlesum\of n =
\frac{11\sqrt{\pi n}}6 -1 +\OO\of{n^{-1/2+\epsilon}}.
\end{equation}

\subsubsection{The asymptotics of $S\of{n,a,b}$
for $a$, $b$ fixed, $n\to\infty$} 
Basically, we mimic the computations
from section \ref{sec:asymSsingle}.
%
In doing so, we are led to 
consider the following function
\begin{equation*}
g\of{n,a,b}\defeq
\sum_{k\geq 1}\sum_{l\geq 1}k^a\,l^b\,\ndivs{\gcd\of{k,l}}\e^{-\pas{k^2+l^2}/n}.
\end{equation*}
Observe again that the terms for $k,l\geq n^{1/2+\epsilon}$
are negligible in this sum.
For obtaining the following formula, we made use of the fact that
$g\of{n,2a,2b}=\OO\of{n^{a+b+1}}$ (which is shown in
appendix~\ref{sec:gnab-asym}):
\begin{multline}
\label{eq:doublesumIng}
\doublesum\of n=
\frac{1}{2}
   \left(-\frac{192}{n^4}+\frac{1632}{n^5}-\frac{8736}{n^6}\right)
   g(n,0,0)
   +\left(\frac{768}
   {n^5}-\frac{8928}{n^6}+\frac{61744}{n^7}\right) g(n,2,0)\\
   +\frac{1}{2}
   \left(-\frac{1152}{n^6}+\frac{18432}{n^7}-\frac{161488}{n^8}\right)
   g(n,2,2)+\left(-\frac{576
   }{n^6}+\frac{10336}{n^7}-\frac{99336}{n^8}\right)
   g(n,4,0)\\
   +\left(\frac{384}
   {n^7}-\frac{10784}{n^8}+\frac{138128}{n^9}\right)
   g(n,4,2)+\frac{1}{2}
   \left(\frac{512}{n^8}-\frac{7872}{n^9}+\frac{58368}{n^{10}}\right)
   g(n,4,4)\\
   +\left(\frac{128}
   {n^7}-\frac{4192}{n^8}+\frac{300624}{5
   n^9}\right)
   g(n,6,0)+\left(-\frac{256
   }{n^8}+\frac{6848}{n^9}-\frac{1517888}{15
   n^{10}}\right)
   g(n,6,2)\\
   +\left(\frac{2432
   }{3 n^{10}}-\frac{62368}{5 n^{11}}\right)
   g(n,6,4)+\frac{1}{2}
   \left(\frac{256}{3 n^{11}}-\frac{416}{15
   n^{12}}\right)
   g(n,6,6)+\left(\frac{544}
   {n^9}-\frac{225488}{15 n^{10}}\right)
   g(n,8,0)\\
   +\left(\frac{3989
   12}{15 n^{11}}-\frac{960}{n^{10}}\right)
   g(n,8,2)+\left(\frac{3173
   6}{15 n^{12}}-\frac{256}{3 n^{11}}\right)
   g(n,8,4)+\frac{1328
   g(n,8,6)}{45
   n^{13}}+\frac{64
   g(n,8,8)}{9
   n^{14}}\\
   +\left(\frac{8672}{5
   n^{11}}-\frac{64}{3 n^{10}}\right)
   g(n,10,0)+\left(\frac{128
   }{3 n^{11}}-\frac{47504}{15 n^{12}}\right)
   g(n,10,2)-\frac{7856
   g(n,10,4)}{45
   n^{13}}-\frac{32
   g(n,10,6)}{3
   n^{14}}\\
   -\frac{456
   g(n,12,0)}{5
   n^{12}}+\frac{2576
   g(n,12,2)}{15
   n^{13}}+\frac{64
   g(n,12,4)}{9
   n^{14}}+\frac{16
   g(n,14,0)}{9
   n^{13}}\\
   -\frac{32
   g(n,14,2)}{9 n^{14}}
   +\OO\of{n^{-5+\epsilon}g\of{n,0,0}}.
\end{multline}

From the results obtained in appendix~\ref{sec:gnab-asym} and
\ref{sec:jacobitheta}, we
easily derive the following asymptotic expansions (using the ``implicit'' definition of
 the numbers $c_{a,b}$ given in \eqref{eq:Zdefinition}):
\begin{align}
g\of{n,0,0}&=
	\frac{\pi ^3 n}{24}
	+\frac{\sqrt{\pi n}}{4} \left(2 
   c_{0,0}- \log
   (n)- \psi\of{\frac{1}{2}}-2
   \gamma  \right)+\OO\of{n^{-m}},
   \notag\\
g\of{n,2a,0}&=
	\frac{2^{-2 a-3} n^{a}  (2 a)!}{ a!}
	\Bigl(
		\frac{\pi ^3n }{3}
		+
   		\sqrt{\pi n} \left(4
   c_{a,0}-\log (n)-\psi\of{a+\frac{1}{2}}-2 \gamma \right)
   \Bigr)+\OO\of{n^{-m}},
   \notag\\
g\of{n,2a,2b} &=
	\frac{2^{-2a-2b-3} n^{a+b}  (2 a)!(2 b)!}{ a!b!}
	\pas{
		\frac{\pi^3 n }{3}
		+
   		4 \sqrt{\pi n}\,
		\frac{
			c_{a,b}\pas{2a+2b}!
		}{
			\pas{a+b}!
		}
	} +\OO\of{n^{-m}}
	\label{eq:asymgnab}
\end{align}
 for all $m\geq 0$. 
Inserting the information from \eqref{eq:asymgnab} in
\eqref{eq:doublesumIng} shows that all the $\log\of n$--terms
cancel, as well as all evaluations of the digamma function $\psi$
(see appendix~\ref{sec:gnb-asym}).
So we obtain 
the second part of the desired result:
\begin{multline}
\label{eq:h2secondpart}
\frac{\pochhammer{n+1}{2}}{12\pas{2n+1}}
\pochhammer{n+1}{3}\doublesum\of n = \\
\sqrt{\pi n} \bigl(
	-\frac{11}{6}
	-2 c_{0,0}
	+8 c_{1,0}
	-9 c_{1,1}
	-9 c_{2,0}
	\\
	+15 c_{2,1}
	+35 c_{2,2}
	+5 c_{3,0}
	-35 c_{3,1}
   \bigr)+\OO\of{n^{-1/2+\epsilon}}.
\end{multline}
Inserting the expressions \eqref{eq:h2firstpart} and
\eqref{eq:h2secondpart} in \eqref{eq:all-2-watermelons-normal} gives the desired result \eqref{eq:2-watermelons-average-asympt}.\hfill\qedsymbol

\begin{appendix}
\section{Background information and relevant results}
\label{sec:appendix}

\subsection{Stirling's approximation applied to quotients of binomial
coefficients}
\label{sec:stirlingsformula}
We have the following asymptotic series for $\log\of{\Gamma\of z}$,
valid for $\absof{\arg z}<\pi-\delta$, $0<\delta<\pi$,
$\absof{z}\to\infty$
(see \cite[equation (3.10.7)]{deBruijn:asymptotic}):
\begin{equation}
\label{eq:Stirlings-formula}
\log\of{\Gamma\of z}\approx
\pas{z-\frac12}\log\of z-z+
\frac{\log\of{2\pi}}{2}+
\sum_{j=1}^\infty z^{1-2j}\frac{B_{2j}}{\pas{2j}\pas{2j-1}},
\end{equation}
where $B_j$ denotes the $j$-th Bernoulli number.

Setting $x=\frac{k-a}{n}$ we thus obtain
\begin{align}
\frac{\binom{2n}{n+a-k}}{\binom{2n}{n}}&=
\exp\Biggl(
	-2n\pas{
		\frac{x^2}{2}
		+\frac{x^4}{12}
		+\frac{x^6}{30}
		+\frac{x^8}{56}+\dots
	}
	+
	\pas{
		\frac{x^2}{2}
		+\frac{x^4}{4}
		+\frac{x^6}{6}
		+\frac{x^8}{8}+\dots
	}\notag\\
	&-\frac1{6 n}
	\pas{
		x^2
   		+x^4+
		x^6+x^8+\dots
	}
	+\frac{1}{n^3}
	\pas{
		\frac{x^2}{30}
		+\frac{x^4}{12}
		+\frac{7x^6}{45}
		+\frac{x^8}{4}+\dots
	}\notag\\
	&-\frac{1}{n^5}
	\pas{
		\frac{x^2}{42}
		+\frac{x^4}{9}
		+\frac{x^6}{3}
		+\frac{11x^8}{14}
		+\dots
	}
	+\frac{1}{n^7}
	\pas{
		\frac{x^2}{30}
		+\frac{x^4}{4}
		+\frac{11 x^6}{10}
		+\frac{143x^8}{40}+\dots
	}\notag\\
	&-\frac{1}{n^9}
	\pas{
   	\frac{5x^2}{66}
	+\frac{5 x^4}{6}
	+\frac{91 x^6}{18}
	+\frac{65 x^8}{3}+\dots
	}+\OO\of{x^2n^{-11}}\label{eq:binkoeffquotapprox}
\Biggr).
\end{align}
Note that $\frac{\binom{2n}{n+a-k}}{\binom{2n}{n}}$ is
zero for $\absof x>1$. The approximation given by
\eqref{eq:binkoeffquotapprox} is very good if, say,
$\absof x \leq \frac12$.

\subsection{Integral representations of the exponential function
	and applications}

\subsubsection{The asymptotics of $g\of{n,b}$ for $b$ fixed, $n\to\infty$}
\label{sec:gnb-asym}
Starting
with the formula
\begin{equation}
\e^{-x}=
\frac{1}{2\pi \i}
\int_{c-\i\infty}^{c-\i\infty}
\Gamma\of zx^{-z}\deeof z \text{ for } c>0, x>1,
\end{equation}
(see \cite[(2.4.1)]{andrews:spec-func})
and using
\begin{equation}
\zeta\of z^2 =\sum_{k\geq 1}\ndivs k k^{-z},
\end{equation}
we obtain
\begin{align*}
g\of{n,b} &=
    \sum_{k\geq 1}\frac{\ndivs{k}}{2\pi \i}\int_{c-\i\infty}^{c-\i\infty}
    n^z\,\Gamma\of z k^{b-2z}\deeof z \\
 &= \frac{1}{2\pi \i}\int_{c-\i\infty}^{c-\i\infty}
    n^z\,\Gamma\of z \zeta\of{2z-b}^2
    \deeof z,
\end{align*}
where $c>\frac{b+1}{2}$. Denote the integrand in the above formula
by $G_1\of{b;z}$.

For any fixed positive number $q$ and $\Re\of s\geq -q$, we have
$\zeta\of s=\OO\of{\absof{s}^{q+1/2}}$ as $s\to\infty$. Since
$n^z\,\Gamma\of z$ becomes small on vertical lines, we can shift
the line of integration to the left as far as we want to, if we
take into account the residues of our integrand $G_1\of{b;z}$.
There is a double pole at
$z=\frac{b+1}{2}$, and possibly some simple poles at
$z=0,-1,-2,\dots$.

For obtaining the residues, we use the power series expansion
\begin{equation*}
\zeta\of s - \frac1{s-1} =
\gamma+\sum_{n=1}^\infty \gamma_n\pas{s-1}^n,
\end{equation*}
where $\gamma$ is Euler's constant and 
$
\gamma_n=\lim_{m\to\infty}\pas{
	\sum_{l=1}^m l^{-1}\pas{\log l}^n-\pas{n+1}^{-1}\pas{\log l}^{n+1}
}
$
(see \cite[1.12, (17)]{erdelyi:high-spec-func}), which
gives the Laurent expansion at $\frac{b+1}2$
\begin{equation*}
\zeta\of{2z-b}^2 = 
	\frac{\gamma}{\pas{z-\frac{b+1}2}^2}+
	\frac{1}{4\pas{z-\frac{b+1}2}}+\cdots.
\end{equation*}
Combining this with the series expansions
\begin{align*}
n^z &= n^{z_0}\pas{
	1+\log\of n\pas{z-z_0}+\dots
}, \\
\Gamma\of{z} &= \Gamma\of{z_0}\pas{
	1+\psi\of{z_0}\pas{z-z_0}+\dots
},
\end{align*}
for $z_0=\frac{b+1}2$,
where $\psi\of z$ is the digamma function (i.e., the logarithmic
derivative of the gamma function $\frac{\Gamma^\prime\of z}{\Gamma\of z}$, see
\cite[1.7]{erdelyi:high-spec-func}), we can easily express the residue
of our integrand $G_1\of{b;z}$ at $z=\frac{b+1}{2}$:
\begin{equation}
\label{eq:residue-g-1}
n^{\frac{b+1}{2}}\Gamma\of{\frac{b+1}{2}}
\pas{
	\frac14\log\of n
	+\frac14\psi\of{\frac{b+1}{2}}
	+\gamma
}.
\end{equation}
For our purposes, we need $\psi\of z$ at
positive integral or half--integral values $z$, which can be derived from
the following information:
\begin{align*}
\psi\of{1} &= -\gamma
	&\text{ (see \cite[section 1.7, equation (4)]{erdelyi:high-spec-func}),} \\
\psi\of{\frac12} &= -\gamma -2\log\of 2
	&\text{ (see \cite[p.~104]{noerlund:differenzen}),} \\
\psi\of{z+n} &= \psi\of z +\sum_{j=0}^{n-1}\frac{1}{z+j}
	&\text{ (see \cite[section 1.7, equation (10)]{erdelyi:high-spec-func}).} \\
\end{align*}
The residue of $G_1\of{b;z}$ at $z=-m$ is
\begin{equation}
\label{eq:residue-g-2}
n^{-m}\frac{(-1)^m}{m!}\zeta\of{-2m-b}^2 =n^{-m}\frac{(-1)^m}{m!} \pas{\frac{B_{2m+b+1}}{\pas{2m+b+1}}}^2,
\end{equation}
where $B_k$ denotes the $k$-th Bernoulli number (see
\cite[1.12, (20)]{erdelyi:high-spec-func}). Note that this
number is non--zero only if $b$ is odd or $m=b=0$.

The sum of \eqref{eq:residue-g-1} and \eqref{eq:residue-g-2} for
all $m\geq 0$ gives an asymptotic series for $g\of{n,b}$.

\subsubsection{The asymptotics of $g\of{n,a,b}$ for $a$, $b$ fixed, $n\to\infty$}
\label{sec:gnab-asym}


In the same manner as in section~\ref{sec:gnb-asym}, we obtain
\begin{align*}
g\of{n,a,b} &=
    \sum_{k,l\geq 1}\frac{\ndivs{\gcd\of{k,l}}}{2\pi \i}\int_{c-\i\infty}^{c-\i\infty}
    n^z\,\Gamma\of z k^a\,l^b \pas{k^2+l^2}^{-z} \deeof z \\
 &= \frac{1}{2\pi \i}\int_{c-\i\infty}^{c-\i\infty}
    n^z\,\Gamma\of z \sum_{k,l\geq 1}\ndivs{\gcd\of{k,l}}
    k^a\,l^b \pas{k^2+l^2}^{-z} \deeof z,
\end{align*}
where $c>\frac{a+b+1}{2}$.

For $k$ and $l$ fixed, set $j=\gcd\of{k,l}$. Then we may write
$k=k_1\,j$
and $l=l_1\,j$ with $\gcd\of{k_1,l_1}=1$. This leads to
\begin{align*}
g\of{n,a,b} &=
	\frac{1}{2\pi \i}
	\int_{c-\i\infty}^{c-\i\infty}
		n^z\,\Gamma\of z
		\underset{\gcd\of{k_1,l_1}=1}{\sum _{j,k_1,l_1\geq1}}
			\ndivs{j}\,j^{a+b-2z}\,
			k_1^{a}l_1^{b}
			\pas{k_1^2+l_1^2}^{-z}\deeof z \\
&=
	\frac{1}{2\pi \i}
	\int_{c-\i\infty}^{c-\i\infty}
		n^z\,\Gamma\of z
		\zeta\of{2z-a-b}^2
		\underset{\gcd\of{k_1,l_1}=1}{\sum _{k_1,l_1\geq1}}
			k_1^{a}l_1^{b}
			\pas{k_1^2+l_1^2}^{-z}\deeof z.
\end{align*}

Now we get rid of the constraint $\gcd\of{k_1,l_1}=1$:
\begin{pro} \label{prop:1}
We have the following identity:
\begin{equation}
\underset{\gcd(k,l)=1}{\sum _{k,l\geq1}}
\frac {k^{a}l^{b}} {({k^2} + {l^2})^z}=
\frac {1} {\zeta(2z-a-b)}{\sum _{k,l\geq1}}
\frac {k^{a}l^{b}} {({k^2} + {l^2})^z}.
\end{equation}
\end{pro}
\begin{proof}
By inclusion--exclusion, the left--hand side equals the sum over
{\it all\/} pairs $(k,l)$ minus the sum over
all pairs $(k,l)$, where some prime number
$p$ divides $\gcd(k,l)$, 
plus the sum over
all pairs $(k,l)$, where the product of two different primes
$p_1p_2$, divides $\gcd(k,l)$, and so on: 
\begin{align*}
\underset{\gcd(k,l)=1}{\sum _{k,l\geq1}}
\frac {k^{a}l^{b}} {({k^2} + {l^2})^z}&=
\sum _{k,l\geq1}
\frac {k^{a}l^{b}} {({k^2} + {l^2})^z}
-
\sum _{p\text{ prime}}\sum _{k,l\geq1}
\frac {(k p)^{a}(l p)^{b}} {({(k p)^2} + {(l p)^2})^z}\\
&\kern3cm
+
\sum _{p_1\neq p_2\text{ prime}}
\sum _{k,l\geq1}
\frac {(k p_1 p_2)^{a}(l p_1 p_2)^{b}} {({(k p_1 p_2)^2} +
{(l p_1 p_2)^2})^z}
-+\cdots\\
&\kern-10pt
=\pas{1-
\sum _{p\text{ prime}}\frac {1} {p^{2z-a-b}}+
\sum _{p_1\neq p_2\text{ prime}}\frac {1} {(p_1 p_2)^{2z-a-b}}-+\cdots
}
\sum _{k,l\geq1}
\frac {k^{a} l^{b}} {({k^2} + {l^2})^z}\\
&\kern-10pt
=
\pas{\prod _{p\text{ prime}} ^{}\of{1-\frac {1} {p^{2z-a-b}}}}
\sum _{k,l\geq1}
\frac {k^{a}l^{b}} {({k^2} + {l^2})^z}.
\end{align*}
The product in the last line is the reciprocal of the
Euler product for $\zeta\of{2z-a-b}$ (\cite[p.~225]{euler:introductio}, which proves
the assertion.
\end{proof}
Thus, we arrive at
\begin{equation}
g\of{n,a,b} =
	\frac{1}{2\pi \i}
	\int_{c-\i\infty}^{c-\i\infty}
		n^z\,\Gamma\of z
		\zeta\of{2z-a-b}
		{\sum _{k,l\geq1}}
			\frac{
				k^{a}l^{b}
			}{
				\pas{k^2+l^2}^{z}
			}
		\,\deeof z.
\end{equation}
Denote the integrand in the above formula
by $G_2\of{a,b;z}$.

Again, we may shift the line of integration to the left as far as we want to, if we
take into account the residues of our integrand $G_2\of{b;z}$.
Computing the poles and residues clearly depends on some
information about the double Dirichlet series involved. This
information will be provided in the next (and final) subsection.

\subsection{The Dirichlet series
$\sum _{k,l\geq1}
\frac {k^{2 a}\,l^{2 b}} {\pas{{k^2} + {l^2}}^z}$}
\label{sec:dirichlet-series}
Note that for our purposes, we only need 
$g\of{n,2a,2b}$, so the series we are interested in is
\begin{equation}
\label{eq:defineZ}
Z\of{a,b;s}\defeq
	\sum_{k,l\geq1}\frac{k^{2a}l^{2b}}{\pas{k^2+l^2}^{s}}.
\end{equation}
Clearly, this is closely related to the following function:
\begin{align}
Z_*\of{a,b;s}
	&\defeq
	\sum_{\pas{k,l}\in\Z^2\setminus\setof{0}}
		\frac{k^{2a}l^{2b}}{\pas{k^2+l^2}^{s}}\notag\\
	&=
	4\cdot Z\of{a,b;s}
	+ 2\cdot\Iverson{b=0}\zeta\of{2s-2a}
	+ 2\cdot\Iverson{a=0}\zeta\of{2s-2b}.
	\label{eq:ZZstar}
\end{align}

Informations on the poles and residues of $Z\of{a,b;s}$
could be directly derived from the work of Pierrette
Cassou--Nogu\`es \cite[p.~41ff]{Cassou-Nogues:1983}, but
there is a simpler way by using the reciprocity law
for Jacobi's theta function.	 This reasoning is a generalization
of Riemann's representation
(see \cite[section 1.12, (16)]{erdelyi:high-spec-func})
of $\zeta\of s$,
\begin{equation}
\pi^{-s/2}\Gamma\of{s/2}\zeta\of s = \frac{1}{s-1}+
	\int_1^\infty\pas{t^{1/2\pas{1-s}}+t^{s/2}}t^{-1}\omega\of t\deeof t,
\end{equation}
where
\begin{equation}
\omega\of t = \sum_{n=1}^\infty\e^{-n^2\pi t} = 1/2\pas{\theta_3\of{0,\i t}-1}.
\end{equation}
Here, $\theta_3$ denotes (one variant of) 
Jacobi's theta function
(see \cite[section 13.19, (8)]{erdelyi:high-spec-func2})
\begin{equation}
\label{eq:define_theta}
\theta_3\of{z,t} =\sum_{n=-\infty}^{\infty}q^{n^2}\e^{2 n \i z},
\end{equation}
where $q=\e^{\i \pi t}$.

\subsubsection{Jacobi's theta function}
\label{sec:jacobitheta}
We rewrite \eqref{eq:define_theta}
by setting $t = y$ and $z =\pi x$ for
$x, y\in\C$; i.e.:
\begin{equation}
\vartheta\of{x,y} =\sum_{n=-\infty}^{\infty}
	\e^{2\pi\i\pas{x n+\frac{y}{2}n^2}}
\end{equation}
This series is absolute convergent for all $x$ and all $y$ with
$\Im\of y >0$. Therefore, for fixed $y$ the function $f: z\mapsto\vartheta\of{z,y}$
is an entire function. We have the properties
\begin{align*}
\vartheta\of{x+1,y} &= \vartheta\of{x,y}, \\
\e^{2\pi\i\pas{x+\frac{y}2}}\vartheta\of{x+y,y} &= \vartheta\of{x,y},
\end{align*}
which in fact determine the theta function up to a multiplicative
factor $c\of y$
(see \cite[section 2.3]{kraetzel:analytische}). Moreover,
we have the following functional equation (see \cite[Theorem 2.12,
equation (2.29)]{kraetzel:analytische}):
\begin{equation*}
\vartheta\of{x,y} =
\vartheta\of{\frac{x}{y},-\frac1y}
\e^{-\pi\i\frac{x^2}{y}}
\sqrt{\frac{\i}{y}}
\end{equation*}
Setting
$
\mytheta\of y \defeq \vartheta\of{0,\i y},
$
we obtain as a special case the following reciprocity law, valid
for all $y$ with $\Re\of y > 0$:
\begin{equation}
\label{eq:my-transformation}
\mytheta\of y =
\sum_{n=-\infty}^{\infty}\e^{-\pi n^2 y} =
\sqrt{\frac{1}{y}}\cdot\mytheta\of{\frac1y}.
\end{equation}

$\mytheta\of y$ is a holomorphic function in the half plane
$\Re\of y>0$, with $\Re\of y = 0$ as essential singular line, see \cite[Satz 2.13]{kraetzel:analytische}.

Interchanging summation, differentiation and integration
in the appropriate places, we obtain
\begin{align}
\frac{\pas{-\pi}^{a+b}}{\pi^s}\Gamma\of{s}Z_*\of{a,b;s}
&= \sum_{\pas{k,l}\in\Z^2\setminus{\setof0}}
		\Gamma\of{s}
		\frac{\pas{-\pi}^{a+b}}{\pi^s}
		\frac{
			k^{2a}l^{2b}
		}{
			\pas{k^2+l^2}^{s}
		} \notag \\
&= \sum_{\pas{k,l}\in\Z^2\setminus{\setof0}}
		\pas{\pi\pas{k^2+l^2}}^{-s}\int_0^\infty
		t^{s-1}\e^{-t}
		\pas{-\pi}^{a+b}k^{2a}l^{2b}
		 \deeof t \notag \\
&=\sum_{\pas{k,l}\in\Z^2\setminus{\setof0}}
		\int_0^\infty
			u^{s-1}
			\pas{-\pi k^2}^a\pas{-\pi l^2}^b
			\e^{-\pi\pas{k^2+l^2}u}
			\deeof u \label{eq:simple-transform} \\
&=\sum_{\pas{k,l}\in\Z^2\setminus{\setof0}}
		\int_0^\infty
			u^{s-1}
			\pas{\differ{\e^{-\pi k^2u}}{u}{a}}
			\pas{\differ{\e^{-\pi l^2u}}{u}{b}}
			\deeof u. \notag
\end{align}
For \eqref{eq:simple-transform}, we used 
$
\Gamma\of s =\int_0^\infty t^{s-1}\e^{-t}\deeof t =
	\alpha^s \int_0^\infty u^{s-1}\e^{-\alpha u}\deeof u
$
with $\alpha=\pas{k^2+l^2}\pi$.

Clearly, $Z_*\of{a,b;s} = Z_*\of{b,a;s}$. So w.l.o.g.\ we may
assume $a\geq b$. We have to distinguish the following two
cases, where we assume $a>0$ and $b\geq0$:
\begin{align}
\frac{1}{\pi^s}\Gamma\of{s}Z_*\of{0,0;s}
&=
\int_0^\infty
	u^{s-1}
	\pas{
		\mytheta\of{u}^2-1
	}\deeof u,\label{eq:case1}\\
\frac{\pas{-\pi}^{a+b}}{\pi^s}\Gamma\of{s}Z_*\of{a,b;s}
&=
\int_0^\infty
	u^{s-1}
	\pas{
		\pas{\differ{\mytheta\of{u}}{u}{a}}
		\pas{\differ{\mytheta\of{u}}{u}{b}}
	}\deeof u.\label{eq:case2}
\end{align}

\paragraph{Case 1}
Considering \eqref{eq:case1}, we basically repeat the reasoning
in \cite[p.\ 203]{kraetzel:analytische}. From
\eqref{eq:my-transformation} we get
$\frac1u\mytheta\of{\frac1u}^2=\mytheta\of{u}^2$.
Assuming $\Re\of s>1$, we compute:
\begin{multline}
\int_0^\infty
	u^{s-1}\pas{\mytheta\of{u}^2-1}
	\deeof u
= 
\int_0^1
	u^{s-1}
	\pas{
		\frac1u\mytheta\of{\frac1u}^2
		-1
	}
	\deeof u
+
\int_1^\infty
	u^{s-1}\pas{\mytheta\of{u}^2-1}
	\deeof u\\
=
\int_0^1
	u^{s-1}
	\pas{
		\frac1u\pas{
			\mytheta\of{\frac1u}^2-1
		}
		+\frac1u
		-1
	}
	\deeof u
+
\int_1^\infty
	u^{s-1}\pas{\mytheta\of{u}^2-1}
	\deeof u\\
= -\frac{1}{s}+\frac{1}{s-1}
	+\int_1^\infty t^{-s}\pas{\mytheta\of{t}^2 -1}\deeof  t
	+\int_1^\infty u^{s-1}\pas{\mytheta\of{u}^2 -1}\deeof  u.
	\label{eq:case1-enum}
\end{multline}
The integrals in the last line converge for all $s\in\C$ and
constitute holomorphic functions, so
$Z_*\of{0,0;s}$ only has a simple pole for $s=1$.

\paragraph{Case 2}
Considering \eqref{eq:case2}, we have to adjust the preceding
method appropriately.
For convenience, set
$
\mytheta_a\of{u}\defeq
		\differ{\mytheta\of{u}}{u}{a}.
$
Then we have:
\begin{equation*}
\int_0^\infty
	u^{s-1}\pas{\mytheta_a\of{u}\mytheta_b\of{u}
	}
	\deeof u
= 
\int_0^1
	u^{s-1}\pas{\mytheta_a\of{u}\mytheta_b\of{u}}\deeof u
+
\int_1^\infty
	u^{s-1}\pas{\mytheta_a\of{u}\mytheta_b\of{u}
	}
	\deeof u.
\end{equation*}
Now use \eqref{eq:my-transformation}in the form
$$\mytheta_a\of{u}\mytheta_b\of{u} =
			\pas{
			\differ{\pas{\sqrt{\frac{1}{u}}\mytheta\of{\frac1u}}}{u}{a}
			}
			\pas{
			\differ{\pas{\sqrt{\frac{1}{u}}\mytheta\of{\frac1u}}}{u}{b}
			}
$$
and combine this with the the 
formula
\begin{equation*}
\differ{\pas{
	\sqrt{\frac1y}\cdot
	f\of{\frac1y}
}}{y}{a} =
\pas{-1}^a
\sum _{k=0}^a
	\frac{
		\fallfac{ 2 a }{ 2 k}
	}{
		4^{k}k!
	}
 \pas{
	\brk{\differ{f}{y}{a-k}}\of{\frac1y}
  y^{k-2a-\frac{1}{2}}
 }
\end{equation*}
to obtain
\begin{multline*}
\mytheta_a\of{u}\mytheta_b\of{u}
=
\frac{
	\pas{-1}^{a+b}(2a)!(2b)!
}{
	4^{a+b}a!b!
}
u^{-a-b-1}
+
\\
\pas{-1}^{a+b}
\sum _{k=0}^a
\sum _{j=0}^b
	\frac{
		\fallfac{ 2 a }{ 2 k}\fallfac{ 2 b }{ 2 j}
	}{
		4^{k+j}k!j!
	}
	\pas{
		\mytheta_{a-k}\of{\frac1u}\mytheta_{b-j}\of{\frac1u}
		-\Iverson{a=k\wedge b=j}
	}
  u^{k+j-2a-2b-1}.
\end{multline*}
Now in the same way as before, this gives
\begin{multline}
\int_0^\infty
	u^{s-1}\pas{\mytheta_a\of{u}\mytheta_b\of{u}
	}
	\deeof u
= \\
\frac{
	\pas{-1}^{a+b}(2a)!(2b)!
}{
	4^{a+b}a!b!\pas{s-a-b-1}
}
+
\int_1^\infty
	u^{s-1}\pas{\mytheta_a\of{u}\mytheta_b\of{u}
	}
	\deeof u
+
\pas{-1}^{a+b}\times\\
\sum _{k=0}^a
\sum _{j=0}^b
	\frac{
			\fallfac{ 2 a }{ 2 k}\fallfac{ 2 b }{ 2 j}
		}{
			4^{k+j}k!j!
		}
\int_1^\infty
		\pas{
			\mytheta_{a-k}\of{u}\mytheta_{b-j}\of{u}
			-\Iverson{a=k\wedge b=j}
	}
  u^{2a+2b-k-j-s}	\deeof u.\label{eq:case2-enum}
\end{multline}
Again, the integrals converge for all $z\in\C$ and constitute
holomorphic functions, whence $Z_*\of{a,b;s}$ has only one
simple pole at $s=a+b+1$.

We summarize all this information in the following proposition.
\begin{pro}
\label{pro:dirichlet-series}
For arbitrary nonnegative integers $a$, $b$,
the series $\sum_{k,l\geq1}\frac{k^{2a}l^{2b}}{\pas{k^2+l^2}^{s}}$
is convergent in the half--plane $\Re\of z > a+b+1$ and defines
a meromorphic function $Z\of{a,b;z}:\C\to\C$ with a simple
pole at $z=a+b+1$, where the residue is
$\frac{
	\pi(2a)!(2b)!
}{
	4^{a+b+1}a!b!\pas{a+b}!
}$.

If $a>0$ and $b>0$, this is the only pole.

If $a=0$ (or $b=0$), there is another simple pole at $z=b+\frac12$
(or $z=a+\frac12$), where the residue is
$-\frac14\pas{\Iverson{b=0}+\Iverson{a=0}}$.

Moreover, we have the following information on special
evaluations of $Z\of{a,b;z}$:
\begin{equation}
\label{eq:special-evaluations}
Z\of{a,b;-n} = \frac18\Iverson{a=b=n=0}\text{ for }n=0,1,2,\dots.
\end{equation}
%
For the absolute term $c_{a,b}$ in the Laurent series expansion
of $Z\of{a,b;z}$ at $z_0=a+b+\frac12$, we have the
following formulas:
\begin{align}
c_{0,0} &= -\gamma -1 +
	\frac12\int_{1}^{\infty}t^{-\frac12}\pas{\mytheta\of{t}^2-1}\deeof t,
	\notag\\
c_{a,0} &= -\frac{\gamma}{2}-\frac{1}{2} +
	\frac{4^{a-1}a!\pas{\pas{-1}^a+1}}{\pas{2a}!}\int_1^\infty
		t^{a-\frac12}\mytheta_a\of t \mytheta\of t\deeof t
	\notag\\
		&+ \sum_{k=1}^a
			\frac{4^{a-k-1}a!}{k!\pas{2a-2k}!}\int_1^\infty
		t^{a-k-\frac12}\pas{\mytheta_{a-k}\of t \mytheta\of t
		-\Iverson{a=k}}\deeof t
		\text{ for } a>0,
		\notag\\
c_{a,b} &=
			\frac{4^{a+b-1}\pas{a+b}!}{\pas{2a+2b}!}
		\Biggl(
		-2\frac{\pas{2a}!\pas{2b}!}{4^{a+b}a!b!}
		+
		\pas{-1}^{a+b}
		\int_1^\infty
			t^{a+b-\frac12}\mytheta_a\of t\mytheta_b\of t \deeof t
		\notag\\
		&+ \sum_{k=0}^a\sum_{j=0}^b
			\frac{\fallfac{2a}{2k}\fallfac{2b}{2j}}{4^{k+j}k!j!}
		\notag\\
		&\times
		\int_1^\infty
			\pas{
				\mytheta_{a-k}\of t\mytheta_{b-j}\of t
				-\Iverson{a=k\wedge b=j}
			}
			t^{a+b-k-j-\frac12}\deeof t
		\Biggr)\text{ for } a\geq b > 0.
		\label{eq:absolute_terms}
\end{align}

\end{pro}
\begin{proof}
This information is extracted straightforwardly from
\eqref{eq:ZZstar} together with
\eqref{eq:case1}, \eqref{eq:case1-enum}
and \eqref{eq:case2}, \eqref{eq:case2-enum}, respectively.

(The evaluation $\zeta\of{-2n}=0$ for nonnegative integers $n$,
which is needed for \eqref{eq:special-evaluations}, can be
found in \cite[1.13, (22)]{erdelyi:high-spec-func}.)
\end{proof}

\end{appendix}

\bibliography{paper}

\end{document}